\newcount\beg
\newbox\tbox
\newbox\aubox
\newbox\adbox
\newbox\mathbox
\newbox\titbox
\newbox\aabox
\newbox\atbox
\newbox\fpbox
\newbox\abbox
\newbox\refbox
\newdimen\refwidth

\font\twbf=cmbx10 at12pt
 at12pt
 at12pt
\font\sc=cmcsc10
\font\tenex=cmr10

\font\ninerm=cmr9 
\font\nineit=cmti9 
\font\ninesy=cmsy9 
\font\ninei=cmmi9 
\font\ninebf=cmbx9 
\font\sevenrm=cmr7  
 
\font\seveni=cmmi7

\font\sevensy=cmsy7 
\font\fivenrm=cmr5  
\font\fiveni=cmmi5  
\font\fivensy=cmsy5

\def\nine{\textfont0=\ninerm \scriptfont0=\sevenrm
           \scriptscriptfont0=\fivenrm
          \textfont1=\ninei \scriptfont1=\seveni
           \scriptscriptfont1=\fiveni
          \textfont2=\ninesy \scriptfont2=\sevensy
           \scriptscriptfont2=\fivensy
          \textfont3=\tenex \scriptfont3=\tenex
           \scriptscriptfont3=\tenex
          \def\rm{\fam0\ninerm}
          \textfont\itfam=\nineit    
          \def\it{\fam\itfam\nineit}
          \textfont\bffam=\ninebf 
          \def\bf{\fam\bffam\ninebf}
          \normalbaselineskip=11pt
          \setbox\strutbox=\hbox{\vrule height8pt depth3pt width0pt}
          \normalbaselines\rm}

\newcommand{\luz}[1]{\luzno#1?}
\newcommand{\luzno}[1]{\ifx#1?\let\next=\relax\yyy
            \else \let\next=\luzno#1\xxx\fi\next}
\renewcommand{\sp}[1]{\def\xxx{\kern1.7pt}\def\yyy{\kern-1.7pt}\luz{#1}}
\newcommand{\spa}[1]{\def\xxx{\kern1pt}\def\yyy{\kern-1pt}\luz{#1}}

\newcommand{\Abbrevauthors}[1]{\setbox\aabox=\hbox{\sevenrm\uppercase{#1}}}
\newcommand{\Abbrevtitle}[1]{\setbox\atbox=\hbox{\sevenrm\uppercase{#1}}}
\newcommand{\Title}[1]{\setbox\tbox=\hbox{\let\\=\cr 
   \baselineskip14pt\vbox{\twbf\tabskip 0pt plus15cc
   \halign to\hsize{\hfil\ignorespaces \uppercase{##}\hfil\cr#1\cr}}}}

\setbox\abbox=\vbox{\vglue18pt}

\newcommand{\Author}[1]{\setbox\aubox=\hbox{\let\\=\cr 
   \nine\baselineskip12pt\vbox{\tabskip 0pt plus15cc
   \halign to\hsize{\hfil\ignorespaces \uppercase{\spa{##}}\hfil\cr#1\cr}}}
   \global\setbox\abbox=\vbox{\unvbox\abbox\box\aubox\vskip8pt}}

\newcommand{\Address}[1]{\setbox\adbox=\hbox{\let\\=\cr 
   \nine\baselineskip12pt\vbox{\it\tabskip 0pt plus15cc
   \halign to\hsize{\hfil\ignorespaces {##}\hfil\cr#1\cr}}}
   \global\setbox\abbox=\vbox{\unvbox\abbox\box\adbox\vskip16pt}}

\newcommand{\Mathclass}[1]{\setbox\mathbox=\hbox{\footnote[]{1991
   {\it Mathematics Subject Classification}\/: #1}}}

\newcommand{\Maketitlebcp}
   {\unhbox\mathbox\vglue5cc\box\tbox\box\abbox\vskip8pt}

\newcommand{\Abstract}[1]{{\nine{\bf Abstract.} #1}}

\newcommand{\Section}[1]{\vskip-\lastskip\vskip12pt plus2pt minus2pt
   {\bf #1}}

\catcode`@=11
\newcommand{\vfootnote}[1]
   {\insert\footins\bgroup\nine\interlinepenalty\interfootnotelinepenalty
   \splittopskip\ht\strutbox\splitmaxdepth\dp\strutbox\floatingpenalty\@MM
   \leftskip\z@skip\rightskip\z@skip\spaceskip\z@skip\xspaceskip\z@skip
   \textindent{#1}\footstrut\futurelet\next\fo@t}
\catcode`@=12

\documentstyle[twoside]{article}

\input amssym.def
\input amssym.tex

\expandafter\chardef\csname pre amssym.def at\endcsname=\the\catcode`\@
\catcode`\@=11

\catcode`\@=\csname pre amssym.def at\endcsname

\Abbrevauthors{Yu. Bespalov and B. Drabant}

\Abbrevtitle{Braided Differential Calculus}

\newtheorem{theorem}{\sc Theorem}
\newtheorem{proposition}[theorem]{\sc Proposition}

\newtheorem{lemma}[theorem]{\sc Lemma}
\newtheorem{definition}{\sc Definition}
\newtheorem{remark}{\sc Remark}
\newtheorem{example}{\sc Example}
\def\abs{\par\vskip 0.3cm\goodbreak\noindent}
\def\DY#1{{\cal DY}\left({#1}\right)}

\def\cO#1{{{\cal O}(#1)}}

\def\hhchh{{{}^H_H{\cal C}^H_H}}
\def\hhdhh{{{}^H_H{\cal D}^H_H}}
\def\d{{\rm d}}
\def\e{{\rm e}}
\def\E{{\bf 1 \!\! {\rm l}}}
\def\i{{\rm i}}
\def\id{{\rm id}}

\def\p{{\rm p}}

\def\m{{\rm m}}
\def\NN{I\!\!N}
\def\overset#1#2{\mathop{\kern0pt #2}\limits ^{#1}}
\def\underset#1#2{\mathop{\kern0pt #2}\limits _{#1}}
\def\lfl{\leaders\hbox to 1em{\hss \hss}\hfill}

\def\nl{\par\noindent}
\def\leftheadline{\noindent{}\rlap{}\hfil\copy\aabox\hfil}
\def\rightheadline{\noindent\hfill\copy\atbox\hfill{}\llap{}}
\def\Mittefrei#1#2{\hbox to \hsize{#1\hss#2}}

\begin{document}

\advance\voffset by -1truecm
\advance\vsize by 1truecm
\setbox\fpbox=\hbox{\sevenrm
************************************************}
\vbox{\sevenrm\baselineskip9pt\hsize\wd\fpbox
\centerline{***********************************************}
\centerline{BANACH CENTER PUBLICATIONS, VOLUME **}
\centerline{INSTITUTE OF MATHEMATICS}
\centerline{POLISH ACADEMY OF SCIENCES}
\centerline{WARSZAWA 1996}\hfill}

\Title{Bicovariant Differential Calculi\\
       and Cross Products\\
       on Braided Hopf Algebras}

\Author{Yuri\ Bespalov}
\Address{National Academy of Sciences\\
Bogolyubov Institute for Theoretical Physics\\
252 143, Kiev-143, Ukraine\\
e-mail: {\ninerm mmtpitp@gluk.apc.org}}

\Author{Bernhard\ Drabant}
\Address{Department of Mathematics, Katholieke Universiteit Leuven\\
         Celestijnenlaan 200B, 3001 Leuven-Heverlee, Belgium\\
         e-mail: {\ninerm bernhard.drabant@wis.kuleuven.ac.be}}

\Maketitlebcp

\Abstract{In a braided monoidal category ${\cal C}$ we consider Hopf
bimodules and crossed modules over a braided Hopf algebra $H$. We show
that both categories are equivalent. It is discussed that the category
of Hopf bimodule bialgebras coincides up to isomorphism with
the category of bialgebra projections over $H$.
Using these results we generalize the Radford-Majid criterion and show
that bialgebra cross products over the Hopf algebra $H$ are precisely
described by $H$-crossed module bialgebras.
In specific braided monoidal abelian categories
we define (bicovariant) braided differential calculi over $H$
and apply the results on Hopf bimodules to construct
a higher order bicovariant differential calculus over $H$
out of any first order bicovariant differential calculus over $H$.
This object is shown to be a bialgebra with universal properties.
}
\abs


\Section{1. Introduction.}
Cross products and cross coproducts of $k$-bialgabras over a field $k$
have been studied in \cite{Rad}.
Let $H$ be a Hopf algebra and $X$ be
an $H$-right module algebra and an
$H$-right comodule coalgebra. Then the conditions are derived
for $H\otimes X$ to be a bialgebra such that the algebra structure is
given by the cross product and the coalgebra structure is given by
the cross coproduct. In this case Radford calls $(H,X)$ an admissible
pair or simply bialgebra cross product over $H$ \cite{Rad}.
The possibility of forming bialgebra cross products over a
quasitriangular Hopf algebra with any of its modules has been recovered in
\cite{Ma2}. In Radford's definition of an admissible pair crossed
modules implicitely appear.
Explicitely crossed modules are defined in \cite{Yet, RT}.
Majid proved that the category of crossed modules
is braided monoidal \cite{Ma1}. He recovered that the 
(co-)modules to be candidates for a bialgebra cross product
are crossed module bialgebras and hence he formulated the Radford
criterion on bialgebra cross products in terms of crossed module
bialgebras.

Bicovariant differential calculi over $k$-Hopf algebras have been
investigated in \cite{Wor}. The so-called bicovariant bimodules \cite{Wor}
or Hopf bimodules are of special importance for the construction of
bicovariant differential calculi. The equivalent description of Hopf
bimodules through crossed modules has been found in \cite{Wor}.

In \cite{Bes, BD1} crossed modules and Hopf (bi-)modules in braided
monoidal categories are constructed.
The equivalence of braided crossed modules and braided Hopf bimodules
has been found in \cite{BD1} as well as the isomorphy of the braided
Hopf bimodule bialgebras and the braided bialgebra projections over a Hopf
algebra $H$ with invertible antipode.
Both theorems can be applied to
generalize the Radford-Majid criterion to braided categories.
In the special case of braided quasitriangular Hopf algebras $\cal H$
\cite{Ma2} this construction works for any $\cal H$-module Hopf algebra 
in a specific admissible category of $\cal H$-modules \cite{Bes, Dra}.
As a result inhomogeneous quantum groups without additional
dilaton generator have been constructed \cite{Dra}.

Another application of the central theorems of \cite{BD1}
in certain braided abelian categories is the construction of
a graded Hopf algebra differential calculus over $H$ which
is in some sense uniquely derived from a given braided bicovariant
first order differential calculus.
This result is discussed in \cite{BD2} and is an
extension of the results of \cite{Wor} to braided categories.
 
\advance\voffset by 1truecm
\advance\vsize by -1truecm
\pagestyle{myheadings}
\markboth{\leftheadline}{\rightheadline}
The paper in hand is concerned with these subjects.
It is mainly based on \cite{BD1, BD2}. In Section 2 we fix our conventions
and give the necessary definitions needed in the following. We
state the central results of \cite{BD1}.
Section 3 is dedicated to cross product constructions. We formulate the
generalized braided Radford-Majid criterion and consider cross products
over braided quasitriangular Hopf algebras or braided quantum groups.
In Section 4 we restrict to
certain braided abelian monoidal categories which we call
$**$-abelian, and define (bicovariant)
differential calculi in a generalized form. Staring from a
bicovariant first order differential calculus over a braided Hopf algebra
$H$ with invertible antipode we deduce a graded Hopf algebra
differential calculus over $H$.
We outline the universality of this graded bialgebra.
For the derivation of these results extensive use of the
results of Section 2 has to be made.


\Section{2. Braided Crossed Modules and Hopf Bimodules.}
For the definition of a monoidal category we refer to \cite{ML1}.
From Mac Lane's coherence theorem \cite{ML2} it is known that
every monoidal category is equivalent to a strict one.
Hence we restrict most of the considerations to strict monoidal
categories and denote by
${\cal C}:=({\cal C},\otimes ,\E, \Psi )$ a
braided monoidal category \cite{FY,JS} where
$\otimes$ is the tensor product (bifunctor), $\E$ is the unit object
and $\Psi$ is the braiding.
We suppose that the reader is familiar with the notion of
(co-)associative (co-)unital (co-)algebras,
(co-)modules and bi-(co-)modules in monoidal categories 
\cite{Swe, Ma2, Ma1}.
In a braided monoidal category the tensor product of two (co-)algebras
is again a (co-)algebra; the multiplication $\m_{U\otimes V}$ and the unit
$\eta_{U\otimes V}$ of two algebras $(U,\m_U,\eta_U)$ and
$(V,\m_V,\eta_V)$ are given through
$\m_{U\otimes V}=
(\m_U\otimes\m_V)\circ(\id_U\otimes\Psi_{U,V}\otimes\id_V)$ and
$\eta_{U\otimes V}=\eta_U\otimes\eta_V$ respectively.
The coalgebra structure of the tensor product of two coalgebras is obtained
in the dual symmetric manner which means that the order of the
composition of morphisms will be reversed, and the multiplication $\m$
will be replaced by the comultiplication
$\Delta$ and the unit $\eta $ by the counit $\varepsilon$, and vice versa.
A bialgebra $(B,\m,\eta,\Delta,\varepsilon)$ in a braided monoidal category
$\cal C$ is an algebra $(B,\m,\eta)$ and a coalgebra
$(B,\Delta,\varepsilon)$ where $\Delta$ and $\varepsilon$ are algebra
morphisms \cite{Ma1}. A Hopf algebra
$(H,\m,\eta,\Delta,\varepsilon,S)$ in $\cal C$ is a bialgebra together with
the antipode $S:H\to H$ such that
$\m\circ(\id_H\otimes S)\circ\Delta=\m\circ(S\otimes \id_H)\circ\Delta=
  \eta\circ\varepsilon$.
Every bialgebra $(B,\m,\eta,\Delta,\varepsilon)$ in $\cal C$ is a 
bi-(co-)module through the
regular action $\m$ and the regular coaction $\Delta$.
The diagonal action of the tensor product of two
right modules $(X,\mu_r^X)$ and $(Y,\mu_r^Y)$ is given by
$\mu_{d,r}^{X\otimes Y}=(\mu_r^X\otimes\mu_r^Y)\circ
(\id_X\otimes\Psi_{Y\,B}\otimes\id_B)\circ(\id_X\otimes\id_Y\otimes\Delta)$.
The right action $\mu_{i,r}$ on $X\otimes Y$ induced by $Y$ is given through
$\mu_{i,r}^{X\otimes Y}=\id_X\otimes\mu_r^Y$.
(Dually) analogue all other types of diagonal (co-)actions and induced
(co-)actions are defined.
We will also use the notion of (co-)adjoint (co-)action.
Let $H$ be a Hopf algebra in $\cal C$ and let $(X,\mu_r,\mu_l)$ be an
$H$-bimodule then $X$ becomes a
right $H$-module through the right adjoint action
${}_{{\rm ad}}\triangleleft:=
\mu_l\circ(\id_H\otimes\mu_r)\circ(\Psi_{X\,H}\otimes \id_H)\circ
   (\id_X\otimes (S\otimes\id_H)\circ\Delta)$.
Similarly the left adjoint action is defined. The coadjoint coactions are
obtained in the dual symmetric manner.
If $A$ is an algebra and $f:H\to A$ is an algebra morphism then the algebra
$A$ becomes an $H$-bimodule $(A,\mu^f_r,\mu^f_l)$
via pullback along $f$, $\mu^f_l=\m_A\circ(f\otimes\id_A)$ and
$\mu^f_r=\m_A\circ(\id_A\otimes f)$.
The corresponding right adjoint action (induced by $f$) will be denoted by
${}_{{\rm ad}_f}\triangleleft$ and the resulting
right $H$-module algebra by $A_f:=(A,\m,\eta,{}_{{\rm ad}_f}\triangleleft)$.
 
We suppose that the categories $\cal C$
admit split idempotents, i.e.
every idempotent in the category $\cal C$, $\e=\e^2:X\to X$ splits
in $\cal C$ in the sense that there exists an object
$X_\e$ and morphisms $\i_\e:X_\e\to X$ and $\p_\e:X\to X_\e$ such that
$\e =\i_\e\circ\p_\e$ and $\id_\e=\p_\e\circ\i_\e$. For example this holds
in every abelian category.
In what follows $B$ and $H$ denote bialgebras and Hopf algebras in
the braided category $\cal C$ respectively. We assume that the antipode
is invertible. The remainder of this section is mostly based on the results
of \cite{Bes, BD1}.

We start with the definition of braided crossed modules.
They are (braided) analogues of the $k$-vector space of invariant
one-forms of a bicovariant bimodule (Hopf bimodule) \cite{Wor} and have
a close connection to smash product constructions (see \cite{Rad} in the
case of bialgebras over the field $k$ and \cite{Bes} for the braided case).

\begin{definition}
{\rm A} right crossed module $(X,\mu_r,\nu_r)$ over the bialgebra $B$
{\rm is
a right $B$-module and a right $B$-comodule obeying the
compatibility relations
\begin{eqnarray}
&&(\id_X\otimes m)\circ(\Psi_{B\,X}\otimes\id_B)\circ
(\id_B\otimes\nu_r\circ\mu_r)\circ
(\Psi_{X\,B}\otimes B)\circ(\id_X\otimes\Delta)
\nonumber\\
&&=(\mu_r\otimes m)\circ(\id_X\otimes\Psi_{B\,B}\otimes\id_B)\circ
(\nu_r\otimes\Delta).
\end{eqnarray}
The category $\DY{\cal C}^B_B$ is the} category of
crossed modules {\rm where the morphisms are both right module and right
comodule morphisms over $B$. In a similar way all other combinations of
crossed modules will be defined.}
\end{definition}

\begin{example}
A Hopf algebra $H$ is a crossed module over itself through
the adjoint action and the regular coaction. Dually analogue
$H$ is a crossed module through
the regular action and the coadjoint coaction.
\end{example}

The (pre-)braided monoidal structure of $\DY{\cal C}^B_B$
is described in the next theorem.

\begin{theorem}\label{braidmon-yd}
For the bialgebra $B$ in $\cal C$ the category
$(\DY{\cal C}^B_B,\otimes,\E)$ is mo\-no\-idal
with unit object $(\E, \varepsilon_B,\eta_B)$.
It is
pre-braided through
\begin{equation}\label{braid-yd}
{}^{{}^{\DY{\cal C}^B_B}}\!\Psi_{X\,Y}:=
(\id_Y\otimes\mu_r^X)\circ(\Psi_{X\,Y}\otimes\id_B)\circ(\id_X\otimes\nu_r^Y)
\end{equation}
where $X,\,Y\in {\rm Ob}(\DY{\cal C}^B_B)$.
If $H$ is a Hopf algebra (with isomorphic antipode) in $\cal C$ then
$\DY{\cal C}^H_H$ is braided, i.e. the inverse of
(\ref{braid-yd}) exists and equals
\begin{eqnarray}
({}^{{}^{\DY{\cal C}^B_B}}\!\Psi_{X\,Y})^{-1}&=&
(\mu_r^X\otimes\id_Y)\circ(\id_Y\otimes(\Psi_{H\,Y})^{-1})
\circ\\
&&((\Psi_{X\,Y})^{-1}\otimes S^{-1})\circ
(\id_Y\otimes(\Psi_{X\,H})^{-1})\circ(\nu_r^Y\otimes\id_X)\nonumber
\end{eqnarray}
\ \lfl$\square$
\end{theorem}

Hopf bimodules in the braided category $\cal C$ are defined as follows.

\begin{definition}
{\rm An object $(X,\mu_r,\mu_l,\nu_r,\nu_l)$ is called a}
$B$-Hopf bimodule in $\cal C$
{\rm if $(X,\mu_r,\mu_l)$ is a $B$-bimodule, and $(X,\nu_r,\nu_l)$ is a
$B$-bicomodule in the category of $B$-bimodules,
where the regular (co-)action on $B$ and the diagonal
(co-)action on tensor products of modules are used.
Hopf bimodules
together with the $B$-bimodule-bicomodule morphisms
form the} category of Hopf bimodules {\rm
which will be denoted by ${}_B^B{\cal C}^B_B$.}
\end{definition}

A consequence of the bimodule property is given by the
following lemma.

\begin{lemma}\label{idem-pi}
Let $(X,\mu_r,\mu_l,\nu_r,\nu_l)$ be an $H$-Hopf
bimodule. Then the morphism ${}_X\Pi:X\to X$ defined through
${}_X\Pi:=\mu_l\circ(S\otimes\id_X)\circ\nu_l$
is an idempotent.
\lfl$\square$
\end{lemma}

We denote by
${}_X\i:{}_HX\to X$, ${}_X\p:X\to {}_HX$ the morphisms which
split the idempotent ${}_X\Pi$, i.e.
${}_X\i\circ {}_X\p={}_X\Pi$ and ${}_X\p\circ {}_X\i={\rm id}_{{}_HX}$.
Then the assignment
${}_H(-):{}^H_H{\cal C}^H_H\to {\DY{\cal C}}^H_H$ which is given through
${}_H(X):= {}_HX$ for an $H$-Hopf bimodule $X$, and through ${}_H(f)=
{}_Y\p\circ f\circ{}_X\i$ for a Hopf bimodule morphism $f:X\to Y$,
defines a functor. The right crossed module structure on ${}_HX$ is given by
$\mu_r^{{}_HX}={}_X\p\circ\mu_r^X\circ({}_X\i\otimes\id_H)$ and
$\nu_r^{{}_HX}=({}_X\p\otimes\id_H)\circ\nu_r^X\circ{}_H\i$.
Conversely if $Y$ is a right crossed module over $H$ then a full inclusion
functor $H\!\ltimes\!(-):
{\DY{\cal C}}^H_H\to {}^H_H{\cal C}^H_H$ is defined by
$H\ltimes (Y) =(H\otimes Y,\mu_{i,l}^{H\otimes Y}, \nu_{i,l}^{H\otimes Y},
\mu_{d,r}^{H\otimes Y}, \nu_{d,r}^{H\otimes Y})$ and
$H\ltimes(f)= \id_H\otimes f$ for any crossed module morphism $f$.
These facts will be used to formulate the following important theorems.

\begin{theorem}\label{Hopf-br}
The category of Hopf bimodules over $H$ is monoidal
with the tensor product given by
$$\otimes_H := \otimes\circ(\id\times {}_H(-)):
{}^H_H{\cal C}_H^H\times{}^H_H{\cal C}_H^H\longrightarrow
{}^H_H{\cal C}_H^H\,.
$$
$H$ is the unit object. It is equipped with the regular $H$-Hopf bimodule
structure.
${}^H_H{\cal C}_H^H$ is braided
with the braiding given on the objects $X,\,Y\in
{\rm Ob}({}^H_H{\cal C}_H^H)$ through
\begin{equation}\label{Hopf-braid}
{}^{{}^{{}^H_H{\cal C}_H^H}}\!\Psi_{X\,Y}
=(\mu_l^Y\otimes{}_X\p\circ\mu_r^X)\circ
 (\id_H\otimes\Psi_{X\,Y}\otimes\id_H)\circ
 (\nu_l^X\otimes\nu_r^Y\circ{}_Y\i)\,.
\end{equation}
\ \lfl$\square$
\end{theorem}

\begin{theorem}\label{yd-hopfbi}
Let $H$ be a Hopf algebra in $\cal C$ with isomorphic antipode. Then the
categories $\DY{\cal C}^H_H$ and ${}^H_H{\cal C}^H_H$ are
braided monoidal equivalent through the functors
$\displaystyle \DY{\cal C}^H_H
{\overset{H\ltimes (-)}{\overrightarrow{\ \ \ \ \ \ \ }}\atop
 \underset{{}_H(-)}{\overleftarrow{\ \ \ \ \ \ \ }}}
{}_H^H{\cal C}^H_H$.\lfl$\square$
\end{theorem}


\Section{3. Braided Cross Products.}
In this section we define braided cross products and generalize
the Radford-Majid criterion to braided categories. We review braided
quasitriangular Hopf algebras and study braided quasitriangular bialgebra
cross products. (See \cite{Bes,BD1,Dra}) 

For the investigation of cross product Hopf algebras,
Hopf algebra projections and Hopf bimodule Hopf algebras
we introduce the relative antipode of a Hopf bimodule.

\begin{definition}\label{rel-ant}
{\rm Let $(X,\mu_r,\mu_l,\nu_r,\nu_l)$ be an $H$-Hopf bimodule.
Then the} relative antipode $S_{X/H}$ of $X$ w.r.t. $H$
{\rm is defined by
\begin{equation}
S_{X/H}:=M_X\circ(S\otimes \id_X\otimes S)\circ N_X
\end{equation}
where $M_X:=\mu_l\circ(\id_H\otimes\mu_r)$ and
$N_X:=(\id_H\otimes\nu_r)\circ\nu_l$.}
\end{definition}

A ``polarized'' form of the anti-(co-)multiplicity of the antipode
      holds,
$S_{X/H}\circ\mu_r=\mu_l\circ\Psi_{X\,H}\circ(S_{X/H}\otimes S)$ and
$S_{X/H}\circ\mu_l=\mu_r\circ\Psi_{H\,X}\circ(S\otimes S_{X/H})$,
and dually analogue for the coactions.
The relative antipode $S_{H/H}$ coincides with the antipode $S$ of $H$.
In the next definition we define braided bialgebra projections in the
sense of \cite{Rad}.

\begin{definition} {\rm Let $B_1$ and $B_2$ be bialgebras
in $\cal C$ and
$B_1\buildrel {\underline\eta}\over\rightarrow
 B_2\buildrel{\underline\varepsilon}\over\rightarrow B_1$
be a pair of bialgebra morphisms such that ${\underline\varepsilon}\circ
{\underline\eta}=\id_{B_1}$. Then $(B_1,B_2,{\underline\eta},
{\underline\varepsilon})$ is called a} bialgebra projection on $B_1$.
{\rm The bialgebra projections on $B_1$ constitute the category
$B_1$-Bialg-$\cal C$. The morphisms obey the relations
$f\circ \underline\eta =\underline\eta$ and
$\underline\varepsilon\circ f=\underline\varepsilon$.}
\end{definition}

One easily verifies that for a bialgebra projection
$(B_1,B_2,\underline\eta,\underline\varepsilon)$
on $B_1$ the object
$\underline B_2=(B_2,\mu_r^{B_2},
\mu_l^{B_2},
\nu_r^{B_2},\nu_l^{B_2})$ is a $B_1$-Hopf bimodule through the (co-)actions
\begin{eqnarray}
&&\mu_l^{B_2}=\m_{B_2}\circ(\underline\eta\otimes\id_{B_2})\,,\quad
  \mu_r^{B_2}=\m_{B_2}\circ(\id_{B_2}\otimes\underline\eta)\,, \\
&&\nu_l^{B_2}=(\underline\varepsilon\otimes\id_{B_2})\circ\Delta_{B_2}\,,\quad
  \nu_r^{B_2}=(\id_{B_2}\otimes\underline\varepsilon)\circ\Delta_{B_2}\,.
\nonumber
\end{eqnarray}
For the formulation of the next theorem observe that in the
category $\hhchh$ the tensor product $\underset H\otimes$ and the
cotensor product $\underset H\square$ exist and coincide up to isomorphism
with $\otimes_H$. The corresponding coequalizer $\lambda$ and equalizer
$\rho$ for two $H$-Hopf bimodules $N$ and $M$ are given respectively by
\begin{eqnarray}\label{tsr-prod}
 &&N\otimes M\overset{\lambda^H_{N,M}}
{\overrightarrow{\ \ \ \ \ \ \ }}N\underset{H}{\otimes}M
\cong N\otimes{}_HM\nonumber\\
&&\lambda^H_{N,M} = (\mu_r^N\otimes\id _{{}_HM})\circ(\id_N\otimes
(\id_H\otimes{}_M\p)\circ\nu_l^M)
\end{eqnarray}
and
\begin{eqnarray}\label{cotsr-prod}
&&N\otimes{}_HM\cong N\underset{H}{\square} M \overset{\rho^H_{N,M}}
{\overrightarrow{\ \ \ \ \ \ \ }}
N\otimes M\nonumber \\
&&\rho^H_{N,M} = (\id_N\otimes\mu_l^M\circ(\id_H\otimes{}_M\i)\circ
 (\nu_r^N\otimes\id_{{}_HM})\,.
\end{eqnarray}
The following theorem states that the categories of Hopf bimodule
bialgebras and bialgebra projections are equal up to isomorphism.

\begin{theorem}\label{hbb-bp}
For any
bialgebra projection $(H,B,\underline\eta_B,\underline\varepsilon_B)$ the
object
$\underline B$ is a bialgebra
$(\underline B,\underline\m_B,\underline\eta_B,\underline\Delta_B,
\underline\varepsilon_B)$ in ${}^H_H{\cal C}^H_H$ where the projection
morphisms
$\underline \eta_B$ and $\underline\varepsilon_B$ are the unit and counit
respectively. Multiplication $\underline\m_B$
and comultiplication $\underline \Delta_B$ are defined through
\begin{equation}\label{bialg-trans1}
\underline\m_B=\m_B\circ(\id_B\otimes{}_H\i)\quad\hbox{and}\quad
\underline\Delta_B=(\id_B\otimes{}_H\p)\circ\Delta\,.
\end{equation}
Conversely every bialgebra $\underline B=(B,\underline\m_B,\underline\eta_B,
\underline\Delta_B,\underline\varepsilon_B)$ in ${}^H_H{\cal C}^H_H$ can be
turned into a bialgebra $B=(B,\m_B,\eta_B,\Delta_B,\varepsilon_B)$ in $\cal C$
where the structure morphisms are given by
\begin{equation}\label{bialg-trans2}
\m_B=\underline\m_B\circ\lambda^H_{B,B}\,,\quad
  \eta_B=\underline\eta_B\circ\eta_H\,,\quad
  \Delta_B=\rho^H_{B,B}\circ\underline\Delta_B\,,\quad
  \varepsilon_B=\varepsilon_H\circ\underline\varepsilon_B\,.
\end{equation}
The pair $(H,B)$ is a bialgebra projection on $H$ via the morphisms
$(\underline\eta_B,\underline\varepsilon_B)$.
There also exists
a correspondence of Hopf algebra structures. For a Hopf algebra
projection $(H_1,H_2,\underline\eta_{H_2},\underline\varepsilon_{H_2})$ on
$H_1$ the antipode of $\underline H_2$ is given by
\begin{equation}\label{anti1}
\underline S_{H_2}=M_{H_2}\circ(\id_{H_1}\otimes S_{H_2}\otimes\id_{H_1})
\circ N_{H_2}
\end{equation}
and for any Hopf algebra
$\underline {H_2}$ in ${}^{H_1}_{H_1}{\cal C}^{H_1}_{H_1}$ the antipode of
${H_2}$ is given by
\begin{equation}\label{anti2}
S_{H_2}=\underline S_{H_2}\circ S_{{H_2}/H_1}=
S_{{H_2}/H_1}\circ \underline S_{H_2}\,.
\end{equation}
In other words the categories $H$-{\rm Bialg}-\,$\cal C$ and
{\rm Bialg}-\,${}^H_H{\cal C}^H_H$
are isomorphic
through the functorial assignment
on the objects described in equations
(\ref{bialg-trans1}) -- (\ref{anti2})
and through the identity on the morphisms.
\nl$\ $\lfl$\square$
\end{theorem}

The equivalence of the braided categories ${}_H^H{\cal C}_H^H$ and
$\DY{\cal C}_H^H$ allows us to modify the previous theorem
in terms of the category $\DY{\cal C}_H^H$.
Then we obtain a formulation of bialgebra cross products in terms
of braided crossed module bialgebras and a braided
generalization of the Radford-Majid criterion
(see \cite{Rad,Ma1} and \cite{Bes,BD1} in the braided case).
This will be outlined in more detail in the following.

The braided cross product and the braided cross coproduct in the category
$\cal C$ can be defined in a formal manner similar to the symmetric
case. For a Hopf algebra $H$ in
$\cal C$ and an algebra $(A,\mu_A)$ in ${\cal C}_H$ the cross product
$H\ltimes_{\!\mu_A} A$ is defined as the universal algebra in
$\cal C$ such that:
\begin{enumerate}
\item
There are algebra morphisms
${\rm j}:H\to H\ltimes_{\!\mu_A} A$ and
${\i}:A\to H \ltimes_{\!\mu_A} A$.
In addition ${\i}$ is algebra morphism in ${\cal C}_H$ where the
module structure on $H\ltimes_{\!\mu_A} A$ is the right adjoint action
induced by the morphism ${\rm j}$, i.e.
${\i}\in{\rm Alg}_H(A,(H\ltimes_{\!\mu_A}A)_{\rm j})$
\item
If $U$ is any algebra in $\cal C$ and $g\in {\rm Alg}(H, U)$,
$f\in {\rm Alg}_H(A,U_g)$, then there exists a unique
algebra morphism $g\ltimes_{\!\mu_A} f: H\ltimes_{\!\mu_A} A\to U$
such that $f=(g\ltimes_{\!\mu_A} f)\circ\i$ and
$g=(g\ltimes_{\!\mu_A} f)\circ{\rm j}$.
\end{enumerate}
The cross coproduct $H\ltimes^{\!\nu_C} C$ is defined in the
dual symmetric manner. Both products are unique up to (co-)algebra
isomorphism. The cross product can be realized on the tensor product
$H\otimes A$ through \cite{Ma2,Dra, BD1}
\begin{eqnarray}
\m_{\ltimes} = (\m_H\otimes\m_A)\circ(\id_H\otimes
	(\id_H\otimes\mu_r^A)\circ(\Psi_{A\,H}\otimes\id_H)\circ
	(\id_A\otimes\Delta_H)\otimes\id_A)\,,&& \nonumber\\
\eta_{\ltimes} = \eta_H\otimes\eta_A\,,\quad
\i=\eta_H\otimes\id_A\,,\quad {\rm j}=\id_H\otimes\eta_A\,.
\hbox to 3cm{\ \hfill}&&
\nonumber
\end{eqnarray}
Then $g\ltimes_{\!\mu_A}f=\m_U\circ(g\otimes f)$ for the corresponding
unique morphism which is induced by $g$ and $f$.

If $X$ is at the same time a right $H$-module algebra and a right
$H$-comodule coalgebra such that the smash product and the smash coproduct
(realized on $H\otimes X$) are compatible in such a way that $H\ltimes X:=
(H\otimes 
X,\m_{\ltimes},\eta_{\ltimes},\Delta_{\ltimes},\varepsilon_{\ltimes})$
is a bialgebra in $\cal C$ then we say, in the sense of \cite{Rad}, that
the pair $(H,X)$ is admissible; $X$ is called
$H$-admissible object in $\cal C$.
The category $H$-cp-$\cal C$ is the category of admissible pairs
$(H,X)$ with bialgebra morphisms $h:H\ltimes X\to H\ltimes Y$ such that
$h\circ{\rm j}_X={\rm j}_Y$ and ${\rm k}_Y\circ h={\rm k}_X$,
where ${\rm k}=\varepsilon\otimes\id_H$.

Without difficulties one proofs that the following relations hold
for an admissible pair $(H,X)$ in $\cal C$.
\begin{eqnarray}\label{adm-rels}
&&\varepsilon_X\circ\m_X =\varepsilon_X\otimes\varepsilon_X\,,\quad
\varepsilon_X\circ\mu_r^X=\varepsilon_X\otimes\varepsilon_H\,,\nonumber\\
&&\Delta_X\circ\eta_X=\eta_X\otimes\eta_X\,,\quad
\nu_r^X\circ\eta_X=\eta_X\otimes\eta_H\,,\\
&&\varepsilon_X \circ \eta_X =\id_{\E}\nonumber
\end{eqnarray}
From (\ref{adm-rels}) it follows that
$(H,H\ltimes X, \id_H\otimes\eta_X,\id_H\otimes\varepsilon_X)$
is a bialgebra projection on $H$.
Using the results of \cite{Bes} and Theorems \ref{yd-hopfbi} and
\ref{hbb-bp} this yields a
description of $H$-admissible objects in the category $\cal C$
in terms of crossed module bialgebras.

\begin{theorem}\label{yd-adm}
Let $H$ be a Hopf algebra in $\cal C$ with isomorphic antipode.
Then the category of admissible pairs $H$-{\rm cp}-$\,{\cal C}$ and
the category of $H$-crossed module bialgebras
{\rm Bialg}-$\DY{\cal C}_H^H$ are isomorphic.
The functorial isomorphism is given through
$\displaystyle \hbox{\rm Bialg-}\DY{\cal C}^H_H
{\overset{H\ltimes (-)}{\overrightarrow{\ \ \ \ \ \ \ \ \ }}\atop
 \underset{{}_H(-)\circ F}{\overleftarrow{\ \ \ \ \ \ \ \ \ }}}
H\hbox{\rm -cp-}\,{\cal C}$ where
$F: H$-{\rm cp}-$\,\cal C\to{\rm Bialg}$-$\hhchh$ is the
restriction of the corresponding
isomorphism of Theorem \ref{hbb-bp}.\lfl$\square$
\end{theorem}

Braided quantum groups (quasitrangular Hopf algebras) 
were introduced in \cite{Majid7} and the basic theory was
developed there.  We use slightly modified definitions
\cite{Bes} which reflect the symmetry between the two coalgebra structures
under consideration.
In what follows we set
$\overline{\cal C}:=({\cal C},\otimes,\E,\Psi^{-1})$, which is
the same monoidal category as $\cal C$ but with inverse braiding.
Then we say that a pair of bialgebras (Hopf
algebras) $(A,\overline A)$ in ${\cal C}\times\overline{\cal C}$
and a convolution invertible bialgebra pairing ${\cal R}:\,\E
\rightarrow\overline A_{\rm op}\otimes A$ in $\cal C$
define a {\it braided qua\-si\-tri\-angu\-lar bialgebra (Hopf algebra)}
$(A,\overline A,{\cal R})$, if
$A$ and $\overline A$ only differ in their comultiplication (and antipode),
i.e.\ $A=(A,\m,\eta,\Delta,\varepsilon (,S))$ and
$\overline A=(A,\m,\eta,\overline\Delta,\varepsilon (,\overline S))$,
and if the identity
$(\Psi_{A,A}\circ\overline\Delta)\cdot{\cal R}={\cal R}\cdot\Delta$
holds, where ``$\cdot$'' is the convolution product and
$\overline A_{\rm op}$ is the opposite bialgebra
w.r.t.\ $\overline A$ in $({\cal C},\otimes,\E,\Psi)$.
Similar to the case of ordinary quantum groups it is shown in \cite{Bes}
that the antipodes of the quasitriangular Hopf algebra
$(A,\overline A,{\cal R})$ are invertible to each
other\footnote{This result does not hold in general in \cite{Ma2}
since there $\overline A_{\rm op}$ is not necessarily supposed to be a
Hopf algebra.}. In
particular $\overline S^{-1}=u\cdot S\cdot u^{-1},$ where
$u:=m\circ(\id_A\otimes S)\circ{\cal R}$. 

For a quasitriangular bialgebra $(A,\overline A,{\cal R})$ we define
the category ${\cal C}_\cO{A,\overline A}$ as the full subcategory of
the category of $A$-right modules ${\cal C}_A$
with objects $(X, \mu_r)$ satisfying the identity
\begin{equation}
(\id_A\otimes\mu_r)\circ(\Psi_{X\,A}\otimes\id_A)\circ(\id_X\otimes\Delta)=
(\id_A\otimes\mu_r)\circ({\Psi^{-1}}_{X\,A}\otimes\id_A)\circ
(\id_X\otimes\overline\Delta)\,.
\end{equation}
This category is identified with a full braided monoidal
subcategory of $\DY{\cal C}^A_A$ \cite{Bes, Majid7};
every module $(X, \mu_r)$ in ${\cal C}_\cO{A,\overline A}$ becomes
a crossed module $(X, \mu_r,\nu_r)$ in $\DY{\cal C}^A_A$ through the
coaction $\nu_r:=(\mu_r\otimes\id_A)\circ(\id_X\otimes{\cal R})$.
This embedding allows us to prove the following generalized
bosonization construction \cite{Bes}.

\begin{theorem}\label{qtha-bos}
Let $(A,\overline A,{\cal R}_A)$ be a quasitriangular bialgebra
in $\cal C$ and $(B,\overline B,{\cal R}_B)$ be a quasitriangular
bialgebra in ${\cal C}_{\cO{A,\overline A}}$ then
$(A\ltimes B,\overline A\ltimes\overline B,{\cal R}_{A\ltimes B})$
is a quasitriangular bialgebra in $\cal C$, where $A\ltimes B$ and
$\overline A\ltimes\overline B$ are the bialgebras obtained through the
functor $A\ltimes(-)$ and $\overline A\ltimes (-)$ in 
$\cal C$ and in $\overline {\cal C}$ respectively. The $R$-matrix of the
product is given by
\begin{equation}
{\cal R}_{A\ltimes B}:=(m_{A\ltimes B}\otimes m_{A\ltimes B})\circ
(\iota_B\otimes(\iota_A\otimes\iota_A)\circ{\cal R}_A\otimes \iota_B)
\circ{\cal R}_B
\end{equation}
where $\iota_A:A\to A\ltimes B$ and $\iota_B:B\to A\ltimes B$
denote the canonical algebra monomorphisms.
Analogous results hold in the Hopf algebra context.
\end{theorem}

Theorem \ref{qtha-bos} allows us to specify
the braided version of the Radford-Majid theorem
to the case where the braided groups are equipped with quasitriangular
structures respected by the projections in the following sense.

\begin{definition}
{\rm Let $(A,\overline A,{\cal R}_A)$ and $(H,\overline H,{\cal R}_H)$
be braided quantum groups in $\cal C$.
Then the pair
$\big(\displaystyle A\mathop{\longrightarrow\atop
\longleftarrow}^{i_A}_{p_A}H\,,\,
\displaystyle \overline A\mathop{\longrightarrow\atop
\longleftarrow}^{i_{\overline A}}_{p_{\overline A}} \overline H\big)$
is called} a quantum group projection {\rm if the following holds.
\begin{enumerate}
\item
Both $(A,H,i_A,p_A)$ and $(\overline A,\overline H,i_A,p_A)$
are bialgebra projections in $\cal C$ and $\overline{\cal C}$
respectively.
\item
For $\underline H\in{\rm Ob}({}^A_A{\cal C}^A_A)$ and
$\underline{\overline H}\in{\rm Ob}
({}^{\overline A}_{\overline A}
{\overline {\cal C}}^{\overline A}_{\overline A})$ according to Theorem
\ref{hbb-bp} the idempotents ${}_{\underline H}\Pi$ and
${}_{\underline{\overline H}}\Pi$ defined in Lemma \ref{idem-pi}
coincide.
\item
$$\!\!\!\!
(\id_H\otimes p_A)\circ{\cal R}_H=(i_A\otimes \id_A)\circ{\cal R}_A
\quad{\rm and}\quad
(p_A\otimes \id_H)\circ{\cal R}_H=(\id_A\otimes i_A)\circ{\cal R}_A \,.
$$
\end{enumerate}}
\end{definition}

Then we obtain a Radford-Majid criterion for braided quantum groups.

\begin{theorem}
Let
$\displaystyle (A,\overline A,{\cal R}_A)
  \mathop{\longrightarrow\atop
\longleftarrow}^{(i_A,i_{\overline A})}_{(p_A,p_{\overline A})}
  (H,\overline H,{\cal R}_H)$ be
a braided quantum group projection in $\cal C$.
Then there exists a quantum group $(B,\overline B,{\cal R}_B)$
in the category
${\cal C}_{\cO{A,\overline A}}$ such that
\begin{equation}
(H,\overline H,{\cal R}_H)\simeq
(A\ltimes B,\overline A\ltimes\overline B,{\cal R}_{A\ltimes B})\,.
\end{equation}
\end{theorem}

\abs

\Section{4. Braided Differential Calculus.}
Differential calculi on quantum groups are constructed in \cite{Wor}
where Hopf bimodules or bicovariant bimodules appear as the basic notion
in Woronowicz's approach to bicovariant differential calculi
on quantum groups. The differential bialgebra structure of the exterior
higher order differential calculus, which has been generated by
a first order bicovariant differential calculus, was found in \cite{Brz}
after the more general investigation of differential bialgebras and quantum
groups had been
worked out in \cite{Mal}. These authors all have in common that they work
over the symmetric tensor category of vector spaces over a field $k$.
Braided differential calculi have been considered in
\cite{Ma4} on the quantum plane -- the fundamental (co-)representation of the
underlying (co-)quasitriangular quantum group.
A successful attempt to construct braided bicovariant differential
calculi is described in \cite{IV} where braided $GL_q(n)$-covariant
differential calculi on the braided matrix algebra $BM_q(n)$
and on the quantum hyperplanes have been found.

In this section we outline the definition and the construction
of braided bicovariant differential calculi \cite{BD2}. We show that
every braided first order bicovariant differential calculus over a
Hopf algebra $H$  in $\cal C$ induces a braided exterior bicovariant
differential calculus over $H$. The universality of this construction
is discussed. For the derivation of these facts we make use of the
results of the previous sections. In what follows we focus
our consideration to
abelian, braided monoidal categories
with a bilinear tensor product, such
that $X$ is flat for every object $X$ in the category.
We will call henceforth categories with these properties
flat abelian categories. The restriction to flat abelian
categories is for convenience. All the results in the sequel can be
derived under weaker conditions. The more general case is studied
in \cite{BD2}.

\begin{remark}\label{rem}
{\rm If the category $\cal C$ of the previous sections is flat abelian then
it is a nontrivial but straightforward computation to verify that
in particular the category $\hhchh$ is also flat abelian.}
\end{remark}

Let us consider any flat abelian category $\cal D$.
Let in the following $I$ be either the set $\{0,1\}$ or $\NN_0$.
Then the $I$-graded category
${\cal D}^I$ is the functor category where $I$ is considered to be discrete,
i.e. the objects in ${\cal D}^I$ are of the form $\hat X=(X_0,X_1,\ldots)$
where $X_j\in {\rm Ob}({\cal D})$ and the morphisms $\hat f:\hat X\to
\hat Y$ are given by $\hat f=(f_0,f_1,\ldots)$ where $f_j:X_j\to Y_j$
is a morphism in $\cal D$ (see e.g.\ \cite{ML1}).
It is not difficult to verify

\begin {proposition}\label{I-grade}
${\cal D}^I$ is a flat abelian category.
In particular the tensor product and the braiding are given respectively by
$(\hat X\otimes\hat Y)_j=\bigoplus_{k+l=j}X_k\otimes Y_l$ and
$(\hat \Psi_{\hat X\,\hat Y})_j=\bigoplus_{k+l=j}(-1)^{kl}\Psi_{X_k\,Y_l}$.
\lfl$\square$
\end{proposition}

The category of complexes over $\cal D$ will be defined in the usual manner
-- see e.g.\ \cite{Hus}.

\begin {definition} {\rm
$(\hat X,\hat \d)$ is called a} complex in ${\cal D}^I$ {\rm if
$\hat X$ is an
object in ${\cal D}^I$ and the differential $\hat \d=(\d_0,\d_1,\ldots)$
is given as a
sequence of morphisms $\d_j:X_j\to X_{j+1}$ in $\cal D$ such that
$\d_{j+1}\circ \d_j=0$.}
Morphisms of complexes $\hat f: (\hat X,\hat \d_X)\to(\hat Y,\hat \d_Y)$
{\rm are morphisms in ${\cal D}^I$ such that $f_{j+1}\circ \d_{X\,j}=
\d_{Y\,j}\circ f_j$. The} category of complexes {\rm is denoted by
${}^c{\cal D}^I$.}
\end{definition}

\begin{proposition}
${}^c{\cal D}^I$ is a flat abelian category.
The tensor product and the braiding are given as in Proposition
\ref{I-grade}. The differential of the tensor product of two complexes
is given through
\begin{equation}\label{tensor-diff}
(\d_{\hat X\otimes\hat Y})_j=\sum_{k+l=j}\left(\d_{X\,k}\otimes
\id_{Y_l}+(-1)^k\,\id_{X_k}\otimes \d_{Y_l}\right)\,.
\end{equation}
$\ $\lfl$\square$
\end{proposition}

Through the canonical mappings $X\mapsto \hat X:=(X,0,0,\ldots)
\mapsto (\hat X, \hat\d:=0)$ one obtains the following
braided monoidal inclusions.

\begin{lemma}\label{inclusion}
$${\cal D}\hookrightarrow{\cal D}^I\hookrightarrow{}^c{\cal D}^I\,.$$
$\ $\lfl$\square$
\end{lemma}

Let $H$ be a Hopf algebra in $\cal D$ (with invertible antipode). Then by
use of Remark \ref{rem} and Lemma \ref{inclusion} we obtain

\begin{proposition}
There are braided monoidal isomorphisms
$${}^H_H({{\cal D}^I})^H_H\cong ({\hhdhh})^I\quad{\rm and}\quad
{}^H_H({{}^c{\cal D}^I})^H_H \cong {}^c({\hhdhh})^I\,.$$
$\ $\lfl$\square$
\end{proposition}

For any object $X$ in $\cal D$ we consider its tensor algebra $T_{\cal D}(X)$
in ${\cal D}^{\NN_0}$ with
$(T_{\cal D}(X))_0=\E_{\cal D}$ and $(T_{\cal D}(X))_j=X^{\otimes\,j}$
the $j$-fold tensor product of $X$. Since $\cal D$ is braided
it is possible for any fixed $j\in I$ to define a section
$S_j\ni\sigma\mapsto\sigma_{\cal D}(X)\in {\rm Aut}(X^{\otimes\,j})$
of the symmetric group $S_j$ to the representation of the braid group $B_j$
on $X^{\otimes\,j}$ generated by the braiding $\Psi_{X\,X}$ \cite{Wor}.
We denote by $S_{(k,l)}$ the shuffle permutations in $S_{k+l}$
which shuffle $k$ numbers in $\{1,2,\ldots,k+l\}$ to the first $k$
places without changing their order and the remaining $l$ numbers
to the places $k+1,\ldots,k+l$ without changing their orders.
Then we define
$A_j(X)=\sum_{\sigma\in S_j}(-1)^{l(\sigma)}\,\sigma_{\cal D}(X)$ and
$A_{(k,l)}(X)=\sum_{\sigma\in S_{(k,l)}}(-1)^{l(\sigma)}\,\sigma_{\cal D}(X)$
where $l(\sigma)$ is the length of the permutation $\sigma$.
Similarly we define
$A^{(k,l)}(X)=\sum_{\sigma\in S^{(k,l)}}(-1)^{l(\sigma)}\,\sigma_{\cal D}(X)$
where $S^{(k,l)}$ are the shuffle permutations inverse
to $S_{(k,l)}$. Using braid algebra one can prove the following proposition.

\begin{proposition}\label{tensor-alg}
On $T_{\cal D}(X)$ a Hopf algebra structure is given by
\begin{eqnarray}\label{hopf-eqns}
\m_{(n,m)} &\cong& \id_{X^{\otimes\,n+m}}: X^{\otimes\,n}\otimes
X^{\otimes\,m}\to X^{\otimes\,n+m}\,,\nonumber\\
\eta_0 &=& \id_k,\ \eta_j=0\ for\ j\ne 0\,,\nonumber\\
\Delta_{(n,m)} &\cong& A_{(n,m)}(X): X^{\otimes\,n+m}\to
X^{\otimes\,n}\otimes
X^{\otimes\,m}\,,\\
\varepsilon_0 &=& \id_k,\ \varepsilon_j=0\ for\ j\ne 0\,,\nonumber\\
S_j &=& (-1)^j\,(\sigma_j^0)_{\cal D}(X)\nonumber
\end{eqnarray}
where $\sigma_j^0=\pmatrix{1 &\ldots &j\cr j &\ldots &1}$.
Because of duality a Hopf structure dual to (\ref{hopf-eqns})
can be established on $T_{\cal D}(X)$. Explicitely it is given by
\begin{eqnarray}\label{co-hopf-eqns}
{\overset\circ\m}_{(n,m)} &\cong& A^{(n,m)}(X):
X^{\otimes\,n}\otimes X^{\otimes\,m}\to X^{\otimes\,n+m}\,,\nonumber\\
{\overset\circ\eta}_0 &=& \id_k,\
{\overset\circ\eta}_j=0\ for\ j\ne 0\,,\nonumber\\
{\overset\circ\Delta}_{(n,m)} &\cong& \id_{X^{\otimes\,n+m}}:
X^{\otimes\,n+m}\to X^{\otimes\,n}\otimes X^{\otimes\,m}\,,\\
{\overset\circ\varepsilon}_0 &=& \id_k,\
{\overset\circ\varepsilon}_j=0\ for\ j\ne 0\,,\nonumber\\
{\overset\circ S}_j &=& (-1)^j\,(\sigma_j^0)_{\cal D}(X)\,.\nonumber
\end{eqnarray}
The dual Hopf algebra according to (\ref{co-hopf-eqns})
will be denoted by ${\overset\circ T}_{\cal D}(X)$.
\lfl$\square$
\end{proposition}

The two Hopf algebras are linked in a nice way via the graded morphism
$(\hat A(X))_j := A_j(X) : X^{\otimes\,j}\to X^{\otimes\, j}$.
This yields the following proposition \cite{BD2}.

\begin{proposition}\label{tensor-forms}
$\hat A(X):T_{\cal D}(X)\to {\overset\circ T}_{\cal D}(X)$ with
$(\hat A(X))_j:=A_j(X)$ is a bialgebra morphism in ${\cal D}^I$.
Hence it follows that $\hat A(X)$ induces a Hopf algebra structure in
${\cal D}^I$ on the object
$$T^\wedge_{\cal D}(X):={\rm coim}(\hat A(X))$$.\lfl$\square$
\end{proposition}

\begin{remark}
{\rm If $\cal D$ is a category of vector spaces then
$T^\wedge_{\cal D}(X)$ is a braided version of the exterior tensor algebra,
$T^\wedge_{\cal D}(X)=T_{\cal D}\big/(\ker\hat A(X))$.}
\end{remark}

From now on we suppose that the category $\cal C$ of the previous sections
is flat abelian. The definition of a differential calculus in $\cal C$
is then given by

\begin{definition} {\rm
A complex $(\hat Y,\hat \d)$ in ${}^c{\cal C}^I$ is called a}
differential calculus {\rm if $(\hat Y,\hat \d)$ is an algebra in
${}^c{\cal C}^I$ and the image of the morphism
$\m_{0,j+1}\circ(\id_{Y_0}\otimes \d_j)$ is $Y_{j+1}$ for all $j\in I$.
Let $\hat {\cal B}$ be a bialgebra in ${\cal C}^I$. Then
$(\hat Y,\hat\d)$ is called} $\hat{\cal B}$-left covariant,
$\hat{\cal B}$-right covariant {\rm or}
$\hat{\cal B}$-bicovariant differential
calculus {\rm if it is a differential calculus in the category
${}^{\hat{\cal B}}({{}^c{\cal C}^I})$, $({{}^c{\cal C}^I})^{\hat{\cal B}}$
or ${}^{\hat{\cal B}}({{}^c{\cal C}^I})^{\hat{\cal B}}$ respectively.
\footnote{According to the general notion these are the
corresponding $\hat{\cal B}$-comodule categories.}}
\end{definition}

If $I=\{0,1\}$ we sometimes use the name first order
differential calculus. If $I=\NN_0$ we speak of higher order differential
calculus \cite{Wor}.
These definitions generalize the notations of \cite{Wor} in the following
sense. If $(\hat Y,\hat\d)$ is a differential calculus and a bialgebra
in ${}^c{\cal C}^I$ then it is straightforward to show that $(\hat Y,\hat\d)$
is a braided $Y_0$-bicovariant differential
calculus
or a braided bicovariant differential calculus over $Y_0$. For $I=\{0,1\}$
this is just a braided version of the definition of bicovariant
differential calculi given in \cite{Wor}.
{\par\vskip 0.3cm\goodbreak}
Now let ${\cal D}=\hhchh$ be the category of $H$-Hopf bimodules
and $((H,X),d)$ be a braided bicovariant
first order differential calculus over the Hopf algebra $H$.
Then it follows rather
immediately that $X$ is in particular an $H$-Hopf bimodule.
Hence we can apply the results of Proposition \ref{tensor-alg}
and \ref{tensor-forms} to derive the Hopf algebra $T^\wedge_{\hhchh}(X)$
in $(\hhchh)^{\NN_0}$. To obtain from $T^\wedge_{\hhchh}(X)$ a Hopf
algebra (projection) in ${\cal C}^{\NN_0}$ we apply the corresponding
functor of Theorem \ref{hbb-bp} and call the resulting object $X^\wedge$.
This object $X^\wedge$ turns out to be the braided generalization of
an algebra of exterior forms over the group as explained in \cite{Wor}
for the case of the symmetric category of vector spaces. To find the
differential structure on $X^\wedge$ is a nontrivial matter \cite{BD2}.
We use ideas of \cite{Wor} and exploit strongly the results and techniques
of \cite{BD1} which have been outlined in the previous sections. We
state the result in the following theorem.

\begin{theorem}
On $X^\wedge$ there exists a unique differential $\hat\d$ such that
\begin{enumerate}
\item
The 0.th component of $\hat\d$ is the differential $\d$ of the bicovariant
first order differential calculus $(H,X,\d)$, i.e.\ $\hat\d_0=\d$.
\item
$(X^\wedge,\hat\d)$ is an $\NN_0$-graded Hopf algebra differential
calculus over $H$ and hence in particular a bicovariant differential
calculus over $H$.
\end{enumerate}
$\ $\lfl $\square$
\end{theorem}

The construction of $X^\wedge$ is universal in the
following sense.
Consider bialgebras $\hat Y$ in ${\cal C}^{\NN_0}$
such that $Y_0$ is a Hopf algebra with invertible antipode and
the components $Y_{n+1}$ are
generated through the multiplication $\m^Y$ by $Y_0$ and $Y_1$,
i.e.\ ${\rm im}(m^Y_{1,n})=Y_{n+1}$ $\forall n>1$. Then
we obtain the next proposition.

\begin{proposition}
Let $\hat Y$ be a bialgebra of the above mentioned form, and let
$(f_0,f_1): (Y_0,Y_1)\to (H,X_1)$ be a bialgebra morphism
in ${\cal C}^{(0,1)}$. Then there exists a unique bialgebra morphism
$\hat f:\hat Y\to X^\wedge$ such that $\hat f_0=f_0$ and $\hat f_1=f_1$.
\phantom{xxxxx}\lfl$\square$
\end{proposition}


\end{document}